\documentclass[12pt]{article}
 \usepackage{latexsym}
 \usepackage{amssymb}
 \usepackage{color}
 \usepackage{graphicx}
\usepackage{epstopdf}
 \newtheorem{Theorem}{Theorem}

\newtheorem{Proposition}{Proposition}

\newcommand{\A}{{\cal A}}

\newcommand{\OO}{{\cal O}}

\newcommand{\vv}{{\bf v}}
\newcommand{\x}{{\bf x}}
\newcommand{\y}{{\bf y}}

\newcommand{\e}{{\bf e}}
\newcommand{\0}{{\bf 0}}

\newcommand{\qed}{\nobreak \ifvmode \relax \else
      \ifdim\lastskip<1.5em \hskip-\lastskip
      \hskip1.5em plus0em minus0.5em \fi \nobreak
      \vrule height0.75em width0.5em depth0.25em\fi}

\def \ep{\hbox{ }\hfill$\Box$}

\begin{document}
\title{Positive Semi-Definiteness of Generalized Anti-Circulant Tensors}

\author{Guoyin Li
\footnote{Department of Applied Mathematics, University of New South
Wales, Sydney 2052, Australia. E-mail: g.li@unsw.edu.au (G. Li).
This author's work was partially supported by Australian Research
Council.} \quad Liqun Qi\footnote{Department of Applied
 Mathematics, The Hong Kong Polytechnic University, Hung Hom, Kowloon, Hong
 Kong.
 E-mail: maqilq@polyu.edu.hk (L. Qi).  This author's work was partially supported by the Hong Kong Research Grant Council (Grant No. PolyU
502510, 502111, 501212 and 501913).}  %\quad Bruce Reznick  \footnote{Department of Mathematics, University of Illinois at Urbana-Champaign, Urbana, IL 61801, USA.  E-mail: reznick@math.uiuc.edu (B. Reznick).}
\quad   Qun Wang \footnote{Department of Applied
 Mathematics, The Hong Kong Polytechnic University, Hung Hom, Kowloon, Hong
 Kong.   Email: wangqun876@gmail.com (Q. Wang).} }

\date{\today} \maketitle

%---------------------------------------------------------------------------------Abstract
\begin{abstract}
\noindent  %\vspace{3mm}
 Anti-circulant tensors have applications in exponential data fitting. They are special Hankel tensors.   In this paper, we extend the definition of anti-circulant tensors to generalized anti-circulant tensors by introducing a circulant index $r$ such that the entries of the generating vector of a Hankel tensor are circulant with module $r$.   In the special case when $r =n$, where $n$ is the dimension of the Hankel tensor, the
 generalized anticirculant tensor reduces to the anti-circulant tensor.   Hence, generalized anti-circulant tensors are still special Hankel tensors.   For the cases that $GCD(m, r) =1$, $GCD(m, r) = 2$ and some other cases, including the matrix case that $m=2$, we give necessary and sufficient conditions for positive semi-definiteness of even order generalized anti-circulant tensors, and show that in these cases, they are sum of squares tensors.   This shows that, in these cases, there are no PNS (positive semidefinite tensors which are not sum of squares) Hankel tensors.

\noindent {\bf Key words:}\hspace{2mm} Anti-circulant tensors, generalized anti-circulant tensor, generating vectors, circulant index, positive semi-definiteness.

\noindent {\bf AMS subject classifications (2010):}\hspace{2mm}
15A18; 15A69
  \vspace{3mm}

\end{abstract}

\newpage
\section{Introduction}
\hspace{4mm}

Anti-circulant tensors were introduced in \cite{DQW}.  They are extensions of anti-circulant matrices in matrix theory \cite{Da, Zh}.   They have applications in exponential data fitting \cite{DQW}.   Anti-circulant tensors are Hankel tensors.   Hankel tensors arise from signal processing and some other applications \cite{BB, PDV, Qi15}.   In this paper, we extend anti-circulant tensors to generalized anti-circulant tensors, which are still Hankel tensors, and present conditions for positive semi-definiteness of generalized anti-circulant tensors.

Let $\vv = (v_0, \cdots, v_{(n-1)m})^\top \in \Re^{(n-1)m+1}$, where $m, n \ge 2$.   An $m$th order $n$ dimensional {\bf Hankel tensor} $\A = (a_{i_1\cdots i_m})$ is defined by
$$a_{i_1\cdots i_m} = v_{i_1+\cdots +i_m-m},$$
for $i_1, \cdots, i_m = 1, \cdots, n$.   If
\begin{equation} \label{e1}
v_i = v_{i+r},
\end{equation}
for $i = 0, \cdots, (n-1)m-r$, where $1 \le r \le n$, then $\A$ is called a {\bf generalized anti-circulant tensor} with {\bf circulant index} $r$.   If $r=n$, then $\A$ is an anti-circulant tensor according to \cite{DQW}.

For $\x = (x_1, \cdots, x_n)^\top \in \Re^n$, $\A$ uniquely define a homogeneous polynomial
\begin{equation} \label{e2}
f(\x) \equiv \A \x^{\otimes m} = \sum_{i_1, \cdots, i_m=1}^n a_{i_1\cdots i_m}x_{i_1}\cdots x_{i_m} = \sum_{i_1, \cdots, i_m=1}^n v_{i_1+\cdots +i_m-m}x_{i_1}\cdots x_{i_m}.
\end{equation}
We call such a polynomial a {\bf Hankel polynomial}.

Suppose that $m=2k$ is even.
In (\ref{e2}), if $f(\x) \ge 0$ for all $\x \in \Re^n$, then $f$ is called a {\bf PSD} (positive semi-definite) {\bf Hankel polynomial} and $\A$ is called a {\bf PSD Hankel tensor} \cite{Qi}.  Denote $\0 = (0, \cdots, 0)^\top \in \Re^n$.   If $f(\x) > 0$ for all $\x \in \Re^n, \x \not = \0$, then $f$ and $\A$ are called {\bf positive definite}.  Clearly, the only odd order PSD Hankel tensor is $\OO$ and there are no odd order positive definite Hankel tensors.  Thus, we may assume that $m$ is even in our discussion.   If $f$ can be decomposed to the sum of squares of polynomials of degree $k$, then $f$ is called an {\bf SOS Hankel polynomial} and $\A$ is called an {\bf SOS Hankel tensor} \cite{HLQS,HLQ,LQX,LQY}.   Clearly, an SOS Hankel tensor is a PSD Hankel tensor but not vice versa.  By \cite{Qi15}, a necessary condition for $\A$ to be PSD is that
\begin{equation} \label{e3}
v_{jm} \ge 0,
\end{equation}
for $j = 0, \cdots, n-1$.
Let $\e_i$ be the $i$th unit vector in $\Re^n$. Substitute them to (\ref{e2}).  Then we get (\ref{e3}) directly.
Vector $\vv$ is called the {\bf generating vector} of $\A$.   It may also generate a $(nk-k+1) \times (nk-k+1)$ Hankel matrix $A = (a_{ij})$ by
$$a_{ij} = v_{i+j-2}$$
for $i, j = 1, \cdots, nk-k+1$.   If the associated Hankel matrix $A$ is PSD, then the Hankel tensor $\A$ is called a {\bf strong Hankel tensor} \cite{Qi15}.   In \cite{LQX}, it was proved that an even order strong Hankel tensor is an SOS Hankel tensor.  Then, a necessary condition for $\A$ to be a strong Hankel tensor is that
\begin{equation} \label{e4}
v_{2j} \ge 0,
\end{equation}
for $j=0, \cdots, (n-1)k$.

Let $GCD(m,r)$ denote the greatest common divisor of the two nonnegative integers $m$ and $r$.

When $r$ is odd, for the case that $m=2k, k \ge 1, GCD(m, r) = 1$ and $n \ge r$, we show that $\A$ is PSD if and only if $v_0 = \cdots = v_{r-1} \ge 0$.   In this case, we show that
$$f(\x) = v_0(x_1+\cdots + x_n)^m,$$
and $\A$ is a strong Hankel tensor.   We show that this result is still true for $r=3, n \ge r$ and $m = 6, 12, 18, 30, 42$.

When $r$ is even, for the case that $m=2k, k \ge 1, GCD(m, r) = 2$ and $n \ge r$, we show that $\A$ is PSD if and only if $v_0 = v_2 = \cdots = v_{r-2}$, $v_1 = v_3 = \cdots = v_{r-1}$, and $v_0 \ge |v_1|$.   In these cases, we may write $v_1 = v_0(2t-1)$, where $t \in [0, 1]$.   We show that
$$f(\x) = tv_0(x_1+\cdots x_n)^m +(1-t)v_0\left(x_1-x_2+x_3-\cdots +(-1)^{n-1}x_n\right)^m,$$
and $\A$ is a strong Hankel tensor.  We show that this result is still true in the case that $m=4$, $r=4$ and $n \ge 4$.

 Note that these two results are true in the matrix case for all $r \ge 1$.   In fact, in the matrix case, for even $r$, we show the result is true as long as $2 \le r \equiv 2p \le 2n-4$.   We believe that our results are new even in the matrix case.

The significance of our results is twofold:  On one hand, they enrich the theory of positive semi-definiteness of even order symmetric tensors \cite{CQ, LL, LWZZL, Qi, Qi15, QS, SQ, ZQZ}.  On the other hand, they also contribute to the study of the Hilbert-Hankel problem.   In 1988, Hilbert \cite{Hi} showed that a PSD symmetric homogeneous polynomial of $n$ variables and degree $m=2k$ is always an SOS polynomial only in the following three cases: 1. $m = 2$; 2. $n=2$; and 3. $m =4, n=3$.   For the other pairs of $m=2k \ge 4$ and $n \ge 3$, there are always PSD polynomials which are not SOS polynomials.  We may call such polynomials PNS (PSD non-SOS) polynomials.   Later, various PNS polynomials are found \cite{Re}.   Since each homogeneous polynomial of $n$ variables and degree $m$ is uniquely corresponding to an $m$th order $n$ dimensional symmetric tensor, the Hilbert discovery may also be stated in the tensor language.   In \cite{LQX}, a question is raised: are all the PSD Hankel tensors SOS tensors?   If the answer to this question is ``yes'', then the problem for determining a given even order Hankel tensor is a PSD tensor or not can be solved by a linear semi-definite programming problem \cite{Las,Lau}, and so, can be solved in polynomial time.   In \cite{LQW}, a conjecture is made that there are no sixth order three dimensional PNS Hankel tensors.   Since generalized anti-circulant tensors are Hankel tensors, our results are related to this Hilbert-Hankel problem.

In Section 2, we prove a theorem for circulant numbers.    This will be useful for our further discussion.

In Sections 3 and 4, we discuss the problem for odd $r$ and even $r$ respectively.  Some final remarks are made in Section 5.

Throughout the paper, we use $\e_i$ to denote the $i$th unit vector in $\Re^n$.

\section{A Theorem on Circulant Numbers}
\hspace{4mm}

We have the following theorem.

\begin{Theorem} \label{t1}
Let $M \ge 1$ and $p \ge 2$.  Suppose that we have a sequence $\{ u_j : j = 0, 1, \cdots \}$, satisfying
$$u_{j+p} = u_j,$$
for $j = 0,  1, \cdots$.  If
\begin{equation} \label{e0}
\sum_{j=0}^M \left({M \atop j}\right)(-1)^ju_{i+j} \ge 0,
\end{equation}
for $i= 0, \cdots, p-1$, or
\begin{equation} \label{e01}
\sum_{j=0}^M \left({M \atop j}\right)(-1)^ju_{i+j} \le 0,
\end{equation}
for $i= 0, \cdots, p-1$, then $u_0 = u_1 = \cdots = u_{p-1}$.
\end{Theorem}
\noindent
{\bf Proof}  We may prove this theorem by induction on $M$.  Obviously, it is true for $M=2$.    Suppose that it is true for $M = 2, \cdots, k$.  We now prove that it is true for $M =k+1$.    Define
$$w_i = \sum_{j=0}^k \left({k \atop j}\right)(-1)^ju_{i+j}$$
for $i = 0, \cdots, p-1$.   Then $w_{i+p} = w_i$ for $i = 0, \cdots, p-1$.     Suppose that (\ref{e0}) holds for $M = k+1$. Note that
\begin{eqnarray*}
w_i-w_{i+1} &= & \left(u_0+  \sum_{j=1}^k \left({k \atop j}\right)(-1)^ju_{i+j}\right)-\left( \sum_{j=0}^{k-1} \left({k \atop j}\right)(-1)^ju_{i+j+1}+ (-1)^{k}u_{i+k+1} \right) \\
& = & u_0+ \left(\sum_{j=1}^k \left(\left({k \atop j}\right)+ \left({k \atop j-1}\right)\right) (-1)^ju_{i+j}\right)+ (-1)^{k+1}u_{i+k+1} \\
&= & u_0+ \left(\sum_{j=1}^k \left({k+1 \atop j}\right) (-1)^ju_{i+j}\right)+ (-1)^{k+1}u_{i+k+1} \\
& = & \sum_{j=0}^{k+1} \left({k+1 \atop j}\right)(-1)^ju_{i+j}.
\end{eqnarray*}
Then (\ref{e0}) is equivalent to
$$w_i - w_{i+1} \ge 0$$
for $i= 0, \cdots, p-1$.   This implies that $w_0 = w_1 \cdots = w_{p-1}$.   Thus, either (\ref{e0}) or (\ref{e01}) holds for $M=k$.  By our induction assumption, we have $u_0 = u_1 = \cdots = u_{p-1}$.   Similarly, if (\ref{e01}) holds for $M = k+1$, we may show that $u_0 = u_1 = \cdots = u_{p-1}$.  This proves the theorem.   \ep

\section{The Case that $r$ is Odd}
\hspace{4mm}

\subsection{The Case that $r=1$}
\hspace{4mm}

This case is trivial.  However,  we present the statement of the result for the reader's convenience here, as it covers the sufficiency part of the results for the cases that $r$ is an odd number with $r \ge 3$.

\begin{Proposition} \label{p0}
Suppose that $\A$ is an $m$th order $n$ dimensional generalized anti-circulant tensor with circulant index $1$, where $m = 2k \ge 2$ and $n \ge 2$.
Then $\A$ is PSD if and only if $v_0 \ge 0$.  In this case, we have
\begin{equation} \label{e13}
f(\x) = v_0(x_1+\cdots+x_n)^m.
\end{equation}
and
\begin{equation} \label{e13.xy}
\y^\top A \y = v_0\left(y_1+\cdots+y_{(nk-k+1)}\right)^2,
\end{equation}
where $A$ is the associated Hankel matrix, which implies that $\A$ is a strong Hankel tensor and hence an SOS Hankel tensor.
\end{Proposition}

The proof is trivial and we omit it here.

\subsection{The Case that $GCD(m, r)=1$}
\hspace{4mm}

We have the following theorem.

\begin{Theorem} \label{t2}
Suppose that $\A$ is an $m$th order $n$ dimensional generalized anti-circulant tensor with $m=2k$, $GCD(m, r) =1$, $1 \le r \le n$ and $k \ge 1$.  Then $\A$ is PSD if and only if $v_0 = \cdots =v_{r-1} \ge 0$.  In this case,  we still have (\ref{e13}) and (\ref{e13.xy}), which implies that $\A$ is a strong Hankel tensor and hence an SOS Hankel tensor.
\end{Theorem}
{\bf Proof}
Suppose that $\A$ is PSD.   Let $\x = \e_q - \e_{q+1}$ for $q=1, \cdots, n$, with $\e_{n+1} \equiv \e_1$.    From $f(\x) \ge 0$, we have
\begin{equation} \label{e13.zw}
\sum_{j=0}^m \left({m \atop j}\right)(-1)^jv_{(q-1)m+j} \ge 0,
\end{equation}
for $q=1, \cdots, n$.    Since $GCD(m, r) = 1$, for each $i= 0, \cdots, r-1$, there is an integer $q$, $1 \le q \le n$ such that $(q-1)m = i$, mod$(r)$.    Then
$v_{i+j} = v_{(q-1)m+j}$ for such $i$, $q$ and $j = 0, \cdots, m$.     Thus, (\ref{e13.zw}) implies that
$$\sum_{j=0}^m \left({m \atop j}\right)(-1)^jv_{i+j} \ge 0,$$
for $i= 0, \cdots, r-1$.   Applying Theorem \ref{t1} with $M = m, u_j = v_j$ and $p = r$, we have $v_0 = \cdots = v_{r-1}$.    By (\ref{e3}), $v_0 \ge 0$.  Thus, we have
$v_0 = \cdots =v_{r-1} \ge 0$.

The ``if'' part follows from Proposition \ref{p0}.
\ep

\subsection{The Case that $GCD(m, r) \not =1$}
\hspace{4mm}

The case that $GCD(m, r) \not =1$ and $r$ is odd includes the case that $r =3, m = 6l$ for $l \ge 1$, the case that $r=5, m = 10l$ for $l \ge 1$, etc.
By \cite{LQW}, Theorem \ref{t2} still holds for the case that $m=6$ and $r=3$.   We may see that Theorem \ref{t2} still holds for more cases that $GCD(m, r) \not =1$ and $r$ is odd.

We now assume that $m=6l$, $r=3$ for $l \ge 1$.

In this case, (\ref{e1}) and (\ref{e2}) have the following form:
\begin{equation} \label{e5}
v_i = v_{i+3}
\end{equation}
for $i = 0, \cdots, (n-1)m-3$, and for $\x = (x_1, \cdots, x_n)^\top \in \Re^n$,
\begin{equation} \label{e6}
f(\x) \equiv \A \x^{\otimes m} = \sum_{i_1, \cdots, i_m=1}^n v_{i_1+\cdots + i_m}x_{i_1}\cdots x_{i_m} = v_0f_0(\x) + v_1f_1(\x)+v_2f_2(\x),
\end{equation}
where
\begin{equation} \label{e7}
f_j(\x) = \sum \left\{ x_{i_1}\cdots x_{i_m} : i_1+\cdots + i_m = j, \ {\rm mod}(3), i_1, \cdots, i_m = 1, \cdots, n \right\},
\end{equation}
for $j= 0, 1, 2$.   We may see that
\begin{equation} \label{e8}
f_0(\x)+ f_1(\x) + f_2(\x) = (x_1+\cdots +x_n)^m.
\end{equation}

Since we are concerned about PSD generalized anti-circulant tensors, we may assume that (\ref{e3}) holds, i.e., $v_0 \ge 0$.

\begin{Proposition} \label{p1}
Suppose that $\A$ is an $m$th order $n$ dimensional generalized anti-circulant tensor with circulant index $3$, where $m = 6, 12, 18, 30, 42$ and $n \ge 3$.  Assume that $v_0 \ge 0$.  If $\A$ is PSD, then
\begin{equation} \label{e9}
v_1 + v_2 = 2v_0.
\end{equation}
\end{Proposition}
\noindent
{\bf Proof}  Suppose that $\A$ is PSD and $v_0 \ge 0$.   Then we have $f(1, -1, 0, \cdots, 0) \ge 0$.   Note that
\begin{eqnarray*}
f_0(1, -1, 0, \cdots, 0) & = & \sum \left\{ x_{i_1}\cdots x_{i_m} : i_1 + \cdots + i_m = 0, \ {\rm mod}(3), \ i_1, \cdots, i_m = 1, 2, \ x_1 = 1,\ x_2 = -1\right\} \\
& = &  \sum \left\{ x_{i_1}\cdots x_{i_m} : {\rm the\ number\ of\ }i_j=1\ {\rm is}\ m-p, \ {\rm the\ number\ of\ }i_j=2\ {\rm is}\ p, \right. \\
& &  \ \ \ \ \ \ \ \ \ \ \ \ \ \ \ \ \ \ \ \ \ p=0, 3, \cdots, m,\ x_1 = 1,\ x_2 = -1\} \\
&= & \sum \left\{ (-1)^p \left({m \atop p}\right) : p=0, 3, \cdots, m \right\}.
\end{eqnarray*}
Similarly, we can prove that
$$f_1(1, -1, 0, \cdots, 0) = \sum \left\{ (-1)^p \left({m \atop p}\right) : p=1, 4, \cdots, m-2 \right\}$$
and
$$f_2(1, -1, 0, \cdots, 0) = \sum \left\{ (-1)^p \left({m \atop p}\right) : p=2, 5, \cdots, m-1 \right\}.$$
By direct calculation, for $m = 6, 18, 30, 42$, we have
\begin{equation} \label{e10}
f_0(1, -1, 0, \cdots, 0) < 0.
\end{equation}
Since $\left({m \atop p}\right) \equiv \left({m \atop m-p}\right)$, we have
\begin{equation} \label{e11}
f_1(1, -1, 0, \cdots, 0) = f_2(1, -1, 0, \cdots, 0).
\end{equation}
By (\ref{e8}),
\begin{equation} \label{e12}
f_0(1, -1, 0, \cdots, 0) + f_1(1, -1, 0, \cdots, 0) + f_2(1, -1, 0, \cdots, 0) = 0.
\end{equation}
By (\ref{e6}), (\ref{e10}-\ref{e12}) and $f(1, -1, 0, \cdots, 0) \ge 0$, we have $v_1+v_2-2v_0 \ge 0$.

On the other hand, for $m = 6, 18, 30, 42$,
\begin{eqnarray*}
&& f_0(1, 1, -2,  0, \cdots, 0) \\
& = & \sum \left\{ x_{i_1}\cdots x_{i_m} : i_1 + \cdots + i_m = 0, \ {\rm mod}(3), \ i_1, \cdots, i_m = 1, 2, 3, \ x_1 = 1,\ x_2 = 1,\ x_3=-2 \right\} \\
& = &  \sum \left\{ x_{i_1}\cdots x_{i_m} : {\rm the\ number\ of\ }i_j=1\ {\rm is}\ m-p-q, \ {\rm the\ number\ of\ }i_j=2\ {\rm is}\ q,  \right. \\
& & \ {\rm the\ number\ of\ }i_j=3\ {\rm is}\ p,\ 2p +q =0, \ {\rm mod}(3),\ 0 \le p, q \le m,\ x_1 = 1,\ x_2 = 1,\ x_3 = -2 \} \\
& = &  \sum_{p=0}^m \sum \left\{ x_{i_1}\cdots x_{i_m} : {\rm the\ number\ of\ }i_j=1\ {\rm is}\ m-p-q, \ {\rm the\ number\ of\ }i_j=2\ {\rm is}\ q,  \right. \\
& & \ {\rm the\ number\ of\ }i_j=3\ {\rm is}\ p,\ 2p +q =0, \ {\rm mod}(3),\ 0 \le q \le m,\ x_1 = 1,\ x_2 = 1,\ x_3 = -2 \} \\
&= & \sum_{p=0}^m \sum \left\{ (-2)^p \left({m \atop p}\right)\left({m-p \atop q}\right) : 2p +q =0, \ {\rm mod}(3),\ 0 \le q \le m \right\}.
\end{eqnarray*}
Similarly, we can prove that
$$f_1(1, 1, -2, 0, \cdots, 0) = \sum_{p=0}^m \sum \left\{ (-2)^p \left({m \atop p}\right)\left({m-p \atop q}\right) : 2p +q =1, \ {\rm mod}(3),\ 0 \le q \le m \right\}$$
and
$$f_2(1, 1, -2, 0, \cdots, 0) = \sum_{p=0}^m \sum \left\{ (-2)^p \left({m \atop p}\right)\left({m-p \atop q}\right) : 2p +q =2, \ {\rm mod}(3),\ 0 \le q \le m \right\}.$$
By direct calculation, for $m = 6, 18, 30, 42$, we have
\begin{equation} \label{e10.a}
f_0(1, 1, -2, 0, \cdots, 0) > 0.
\end{equation}
Note that $\left({m-p \atop q}\right) \equiv \left({m-p \atop m-p-q}\right)$.  Also, $2p+q = 2$, mod$(3)$ is equivalent to $2p+m-p-q = 1$, mod$(3)$.   Thus,
\begin{equation} \label{e11.a}
f_1(1, 1, -2, 0, \cdots, 0) = f_2(1, 1, -2, 0, \cdots, 0).
\end{equation}
By (\ref{e8}),
\begin{equation} \label{e12.a}
f_0(1, 1, -2, 0, \cdots, 0) + f_1(1, 1, -2, 0, \cdots, 0) + f_2(1, 1, -2, 0, \cdots, 0) = 0.
\end{equation}

By $f(1, 1, -2, 0, \cdots, 0) \ge 0$, (\ref{e6}) and (\ref{e10.a}-\ref{e12.a}) we can derive $v_1+v_2-2v_0 \le 0$.   This proves (\ref{e9}).

For the case that $m = 12$, (\ref{e11}) still holds.   By direct computation, we have
\begin{equation} \label{e10ab}
f_0(1, -1, 0, 0, \cdots, 0) > 0.
\end{equation}
By (\ref{e6}), (\ref{e11}-\ref{e12}) and (\ref{e10ab}), we have $v_1+v_2-2v_0 \le 0$.
On the other hand, consider $f(1, -3, 2, 0, \cdots, 0)$ and $f(1, 2, -3, 0, \cdots, 0)$.   By direct computation, we have
\begin{equation} \label{e10cd}
f_0(1, -3, 2, 0, \cdots, 0) < 0.
\end{equation}
We have that
\begin{eqnarray*}
&& f_0(1, -3, 2,  0, \cdots, 0) \\
&= & \sum_{p=0}^m \sum \left\{ 2^p(-3)^q \left({m \atop p}\right)\left({m-p \atop q}\right) : 2p +q =0, \ {\rm mod}(3),\ 0 \le q \le m \right\}
\end{eqnarray*}
and
\begin{eqnarray*}
&& f_0(1, 2, -3,  0, \cdots, 0) \\
&= & \sum_{p=0}^m \sum \left\{ 2^q(-3)^p \left({m \atop p}\right)\left({m-p \atop q}\right) : 2p +q =0, \ {\rm mod}(3),\ 0 \le q \le m \right\}
\end{eqnarray*}
We may see that $2p+q =0$, mod$(3)$ if and only if $p+2q =0$, mod$(3)$.  Thus,
\begin{equation} \label{e10ef}
f_0(1, 2, -3, 0, \cdots, 0) = f_0(1, -3, 2, 0, \cdots, 0) < 0.
\end{equation}

Similarly, we may show that
\begin{equation} \label{e10gh}
f_1(1, 2, -3, 0, \cdots, 0)-f_2(1, 2, -3, 0, \cdots, 0) = f_2(1, -3, 2, 0, \cdots, 0)-f_1(1, -3, 2, 0, \cdots, 0).
\end{equation}
By $f(1, -3, 2, 0, \cdots, 0)+f(1, 2, -3, 0, \cdots, 0) \ge 0$, (\ref{e6}) and (\ref{e10ef}-\ref{e10gh}) we can derive $v_1+v_2-2v_0 \ge 0$.   This proves that (\ref{e9}) still holds for $m=12$.
 \ep

We now have the following theorem.

\begin{Theorem} \label{t3}
Suppose that $\A$ is an $m$th order $n$ dimensional generalized anti-circulant tensor with $m=6, 12, 18, 30, 42$, $r=3$ and $n \ge r$.  Then $\A$ is PSD if and only if $v_0 =v_1=v_2 \ge 0$.  In this case,  we still have (\ref{e13}) and (\ref{e13.xy}), which implies that $\A$ is a strong Hankel tensor and hence an SOS Hankel tensor.
\end{Theorem}
\noindent
{\bf Proof}  Suppose that $\A$ is PSD.  Then $v_0 \ge 0$.  Without loss of generality, assume that $v_0 > 0$.    By Proposition \ref{p1}, $v_1 + v_2 = 2v_0$.   Let $v_1 = v_0(1+\alpha)$.   Then $v_2 = v_0(1-\alpha)$ and
\begin{equation} \label{e12.b}
f(\x) = v_0(x_1+\cdots +x_n)^m +v_0\alpha(f_1(\x)-f_2(\x)),
\end{equation}
 where $f_1$ and $f_2$ are defined as in (\ref{e7}).   We may see that
$$f_1(1, 2, -3, 0, \cdots, 0)-f_2(1, 2, -3, 0, \cdots, 0) = f_2(1, -3, 2, 0, \cdots, 0)-f_1(1, -3, 2, 0, \cdots, 0)  \not = 0.$$
Then from this, (\ref{e12.b}),  $f(1, 2, -3, 0, \cdots, 0) \ge 0$ and $f(1, -3, 2, 0, \cdots, 0) \ge 0$, we have $\alpha = 0$.   This proves that $v_0=v_1=v_2 \ge 0$. The remaining conclusions now follow from Proposition \ref{p0}.   \ep

In the proof of Proposition \ref{p1}, we use direct calculation to show (\ref{e10}), (\ref{e10.a}), (\ref{e10ab}) and (\ref{e10cd}).  Are (\ref{e10}) and (\ref{e10.a}) still true for $m = 12l+6$ with $l \ge 4$?   Are (\ref{e10ab}) and  (\ref{e10cd}) still true for $m = 12l$ with $l \ge 2$?    How can we prove these by some analytical technique?  The case that $m=2k, r=2p+1, GCD(m, r) \not = 1$ for $k \ge 2$ and $p \ge 2$ also remains unknown.    These are some further research topics.

\section{The Case that $r$ is Even}
\hspace{4mm}

\subsection{The Case that $r=2$}
\hspace{4mm}

We see that the results in \cite{LQW} for $m=6$ and $n=3$ can be extended to this case.

In this case, (\ref{e1}) and (\ref{e2}) have the following form:
\begin{equation} \label{e3.1}
v_i = v_{i+2}
\end{equation}
for $i = 0, \cdots, (n-1)m-2$, and for $\x = (x_1, \cdots, x_n)^\top \in \Re^n$,
\begin{equation} \label{e3.2}
f(\x) \equiv \A \x^{\otimes m} = \sum_{i_1, \cdots, i_m=1}^n v_{i_1+\cdots + i_m-m}x_{i_1}\cdots x_{i_m} = v_0f_0(\x) + v_1f_1(\x),
\end{equation}
where
\begin{equation} \label{e3.3}
f_j(\x) = \sum \left\{ x_{i_1}\cdots x_{i_m} : i_1+\cdots + i_m = j, \ {\rm mod}(2), i_1, \cdots, i_m = 1, \cdots, n \right\},
\end{equation}
for $j= 0, 1$.   We may see that
\begin{equation} \label{e3.4}
f_0(\x)+ f_1(\x) = (x_1+\cdots +x_n)^m.
\end{equation}

Since we are concerned about PSD generalized anti-circulant tensors, we may assume that (\ref{e3}) holds, i.e., $v_0 \ge 0$.

\begin{Theorem} \label{t4}
Suppose that $\A$ is an $m$th order $n$ dimensional generalized anti-circulant tensor with circulant index $r=2$, where $m=2k \ge 2$ and $n \ge 2$.   Then $\A$ is PSD if and only if $|v_1| \le v_0$.   In these cases, we may write $v_1 = v_0(2t-1)$, where $t \in [0, 1]$.   We have that
$$f(\x) = tv_0(x_1+\cdots x_n)^m +(1-t)v_0\left(x_1-x_2+x_3-\cdots +(-1)^{n-1}x_n\right)^m,$$
and $\A$ is a strong Hankel tensor.
\end{Theorem}
\noindent
{\bf Proof}  Suppose that $\A$ is PSD and $v_0 \ge 0$.   Then we have $f(1, 1, 0, \cdots, 0) \ge 0$.   Note that
\begin{eqnarray*}
f_0(1, 1, 0, \cdots, 0) & = & \sum \left\{ x_{i_1}\cdots x_{i_m} : i_1 + \cdots + i_m = 0, \ {\rm mod}(2), \ i_1, \cdots, i_m = 1, 2, \ x_1 = 1,\ x_2 = 1\right\} \\
& = &  \sum \left\{ x_{i_1}\cdots x_{i_m} : {\rm the\ number\ of\ }i_j=1\ {\rm is}\ m-p, \ {\rm the\ number\ of\ }i_j=2\ {\rm is}\ p, \right. \\
& &  p=0, 2, \cdots, m,\ x_1 = 1,\ x_2 = 1\} \\
&= & \sum \left\{ \left({m \atop p}\right) : p=0, 2, \cdots, m \right\}.
\end{eqnarray*}
Similarly, we can prove that
$$f_1(1, 1, 0, \cdots, 0) = \sum \left\{ \left({m \atop p}\right) : p=1, 3, \cdots, m-1 \right\},$$
$$f_0(1, -1, 0, \cdots, 0) = \sum \left\{ \left({m \atop p}\right) : p=0, 2, \cdots, m \right\}$$
and
$$f_1(1, -1, 0, \cdots, 0) = \sum \left\{ (-1)^p \left({m \atop p}\right) : p=1, 3, \cdots, m-1 \right\}.$$
We may see that
\begin{eqnarray} \label{e3.5}
  &&f_0(1, -1, 0, \cdots, 0)+f_1(1, -1, 0, \cdots, 0)=0,\nonumber\\
  &&f_0(1, 1, 0, \cdots, 0) = f_0(1, -1, 0, \cdots, 0) > 0,\nonumber\\
  &&f_1(1, 1, 0, \cdots, 0) = -f_1(1, -1, 0, \cdots, 0) > 0.
\end{eqnarray}
By (\ref{e3.2}), we have $v_0+v_1 \ge 0$.  From $f(1, -1, 0, \cdots, 0) \ge 0$ and (\ref{e3.2}), we have $v_0-v_1 \ge 0$.   This implies that $v_0 \ge |v_1|$.
On the other hand, suppose that $v_0 \ge |v_1|$.   We may write $v_1 = v_0(2t-1)$, where $t \in [0, 1]$. Write $f(\x)=v_0 f_0(\x)+ v_1 f_1(\x)$ such that $f_0(\x)+f_1(\x)=(x_1+\cdots +x_n)^m$ and $f_0(\x)-f_1(\x)=(x_1-x_2+x_3-\cdots+(-1)^{n-1}x_n)^m$ for all $\x=(x_1,\cdots,x_n)^\top \in \Re^{n}$. It then follows
from (\ref{e3.2})  that
\begin{eqnarray*}
f(\x) = v_0 f_0(\x)+v_1 f_1(\x)
& = & v_0 f_0(\x)+ (2t-1) v_0 f_1(\x) \\
& = & tv_0 (f_0(\x)+f_1(\x))+(1-t) v_0 (f_0(\x)-f_1(\x)) \\
& = & tv_0(x_1+\cdots +x_n)^m + (1-t)v_0(x_1-x_2+x_3-\cdots+(-1)^{n-1}x_n)^m.
\end{eqnarray*}
Similarly, we have
$$g(\y) = \y^\top A \y = tv_0(y_1+\cdots+y_{nk-k+1})^2 + (1-t)v_0(y_1-y_2+y_3-\cdots+(-1)^{nk-k}y_{nk-k+1})^2,$$
where $A$ is the associated Hankel matrix of $\A$.
The conclusions now follow from the definitions of PSD, SOS and strong Hankel tensors.
\ep

\subsection{The Case that $GCD(m, r)=2$}
\hspace{4mm}

In this section, we allow $r \le 2n-4$ instead of $r \le n$, and still call such a tensor a generalized anti-circulant tensor.     We have the following theorem.

\begin{Theorem} \label{t5}
Suppose that $\A$ is an $m$th order $n$ dimensional generalized anti-circulant tensor with $m = 2k, k \ge 1$, $4 \le r=2p \le 2n-4$. If $GCD(m, r) = 2$, then $\A$ is PSD if and only if $v_0 = v_2 = \cdots = v_{r-2}$, $v_1 = v_3 = \cdots =v_{r-1}$ and $v_0 \ge |v_1|$.  In this case, we may write $v_1 = v_0(2t-1)$, where $t \in [0, 1]$.  Then we have
$$f(\x) =  tv_0(x_1+\cdots x_n)^m +(1-t)v_0\left(x_1-x_2+x_3-\cdots +(-1)^{n-1}x_n\right)^m.$$
This implies that $\A$ is PSD if only if it is SOS.  Furthermore, in this case, $\A$ is a strong Hankel tensor.
\end{Theorem}
\noindent
{\bf Proof}   Suppose that $\A$ is PSD.  Let $\x = \e_q - \e_{q+2}$ for $q=1, \cdots, n$, with $\e_{n+1} \equiv \e_1$ and $\e_{n+2} \equiv \e_2$.  By $f(\x) \ge 0$, we have that
\begin{equation} \label{e13.uv}
\sum_{j=0}^m \left({m \atop j}\right)(-1)^jv_{(q-1)m+2j} \ge 0,
\end{equation}
for $q= 1, \cdots, n$.   Since $GCD(m, r) = 2$, for each $i= 0, \cdots, p-1$, there is an integer $q$, $1 \le q \le n$ such that $(q-1)m = 2i$, mod$(r)$.    Then
$v_{2(i+j)} = v_{(q-1)m+2j}$ for such $i$, $q$ and $j = 0, \cdots, m$.     Thus, (\ref{e13.uv}) implies that
$$\sum_{j=0}^m \left({m \atop j}\right)(-1)^jv_{2(i+j)} \ge 0,$$
for $i= 0, \cdots, p-1$.   Applying Theorem \ref{t1} with $M = m$ and $u_j = v_{2j}$, we have $v_0 = v_2 = \cdots = v_{r-2}$.

Let $\x = \alpha \e_{q-1} + \e_q - \alpha \e_{q+1}$ for $q = 1, \cdots, n$ with $\e_0 \equiv \e_n$.    Since $v_0 = v_2 = \cdots = v_{r-2}$,
in the expression of $f(\x)$, the coefficient for power $\alpha^m$ is zero.   Hence, the highest power of $\alpha$ in $f(\x)$ is the term for power $\alpha^{m-1}$, which is
$$m \alpha^{m-1}\left(\sum_{j=0}^{m-1}\left({m-1 \atop j}\right)(-1)^jv_{mq-2m+1+2j}\right).$$
From $f(\x) \ge 0$, letting $\alpha \to \infty$, we have
\begin{equation} \label{e13.st}
\sum_{j=0}^{m-1}\left({m-1 \atop j}\right)(-1)^jv_{mq-2m+1+2j} \ge 0,
\end{equation}
for $q = 1, \cdots, n$.    Since $GCD(m, r) = 2$, as in the first part of this proof, (\ref{e13.st}) implies that
$$\sum_{j=0}^{m-1}\left({m-1 \atop j}\right)(-1)^jv_{2i+1+2j} \ge 0,$$
for $i= 0, \cdots, p-1$.  Applying Theorem \ref{t1} with $M = m-1$ and $u_j = v_{2j+1}$, we have $v_1 = v_3 = \cdots =v_{r-1}$.
The remaining conclusions now follow from Theorem \ref{t4}.  \ep

\subsection{The Case that $GCD(m, r) = 2l$ for $l \ge 2$}
\hspace{4mm}

In this case, we have the following theorem for $m=4$, $n \ge r=4$.

\begin{Theorem} \label{t6}
Suppose that $\A$ is a fourth order $n$ dimensional generalized anti-circulant tensor with circulant index $r=4$, where $n \ge 4$.   Then $\A$ is PSD if and only if $v_0 = v_2, v_1= v_3$ and $|v_1| \le v_0$.   In these cases, we may write $v_1 = v_0(2t-1)$, where $t \in [0, 1]$.   We have that
$$f(\x) = tv_0(x_1+\cdots x_n)^4 +(1-t)v_0\left(x_1-x_2+x_3-\cdots +(-1)^{n-1}x_n\right)^4,$$
and $\A$ is a strong Hankel tensor.
\end{Theorem}
\noindent
{\bf Proof}
In this case, (\ref{e1}) and (\ref{e2}) have the following form:
\begin{equation} \label{e25}
v_i = v_{i+4}
\end{equation}
for $i=0,\cdots,4n-8$. From (\ref{e2}), for $\x=(x_1, \cdots ,x_n)^\top \in \Re^n$, we have
\begin{equation} \label{e26}
f(\x) \equiv \A \x^{\otimes 4} = \sum_{i_1, \cdots, i_4=1}^n v_{i_1+\cdots + i_4}x_{i_1}\cdots x_{i_4} = v_0f_0(\x) + v_1f_1(\x)+v_2f_2(\x)+v_3f_3(\x),
\end{equation}
where
\begin{equation} \label{e26.1}
f_j(\x) = \sum \left\{ x_{i_1}\cdots x_{i_4} : i_1+\cdots + i_4 = j, \ {\rm mod}(4), i_1, \cdots, i_4 = 1, \cdots, n \right\},
\end{equation}
for $j= 0, 1, 2, 3$.
Furthermore, we have
\begin{equation} \label{e36}
f_0(x_1, x_2, x_3, x_4, 0, \cdots 0) = x_1^4 + x_2^4+ x_3^4+ x_4^4 + 6(x_1^2x_3^2+ x_2^2x_4^2)+ 12(x_1^2x_2x_4 + x_1x_2^2x_3+ x_2x_3^2x_4 + x_1x_3x_4^2),
\end{equation}
\begin{equation} \label{e37}
f_1(x_1, x_2, x_3, x_4, 0, \cdots 0) = 4(x_1^3x_2+ x_2^3x_3+ x_3^3x_4+ x_1x_4^3)+ 12(x_1^2x_3x_4+ x_1x_2^2x_4 + x_1x_2x_3^2+ x_2x_3x_4^2),
\end{equation}
\begin{equation} \label{e38}
f_2(x_1, x_2, x_3, x_4, 0, \cdots 0) = 4(x_1^3x_3 + x_1x_3^3+ x_2^3x_4+ x_2x_4^3)+ 6(x_1^2x_2^2+ x_2^2x_3^2+ x_3^2x_4^2+ x_1^2x_4^2)+ 24x_1x_2x_3x_4,
\end{equation}
\begin{equation} \label{e39}
f_3(x_1, x_2, x_3, x_4, 0, \cdots 0) = 4(x_1x_2^3+ x_2x_3^3+ x_3x_4^3+ x_1^3x_4)+ 12(x_1^2x_2x_3+ x_2^2x_3x_4+ x_1x_3^2x_4+x_1x_2x_4^2).
\end{equation}

Suppose now that $\A$ is PSD.   From (\ref{e37}) and (\ref{e39}), we see that
$$f_1(1, 0, -1, 0, \cdots, 0) = f_3(1, 0, -1, 0, \cdots, 0) = 0.$$
From (\ref{e36}) and (\ref{e38}), we have
$$f_0(1, 0, -1, 0, \cdots, 0) = -f_2(1, 0, -1, 0, \cdots, 0) > 0.$$
Then by $f(1, 0, -1, 0, \cdots, 0) \ge 0$, we have $v_0 \ge v_2$.

Similarly, from (\ref{e36}-\ref{e39}), we have
$$f_1(1, -1, -1, 1, 0, \cdots, 0) = f_3(1, -1, -1, 1, 0, 0, \cdots, 0) = 0$$
and
$$f_0(1, -1, -1, 1, 0, \cdots, 0) = -f_2(1, -1, -1, 1, 0, \cdots, 0) < 0.$$
Then by $f(1, -1, -1, 1, 0, \cdots, 0) \ge 0$, we have $v_0 \le v_2$.   Thus, we derive that $v_0=v_2$.

From $f(\alpha, 1, -\alpha, 0, \cdots, 0) \ge 0$,
$f(\alpha, -1, -\alpha, 0, \cdots, 0) \ge 0$ and (\ref{e26.1}), we derive that $v_0 \ge \phi(\alpha)|v_3 -v_1|$, where $\phi(\alpha) \to \infty$ if $\alpha \to \infty$.   Letting $\alpha$ tend to $\infty$, we have $v_1 = v_3$.   The remaining conclusions now follow from Theorem \ref{t4}.  \ep

An interesting question would be: can we extend this theorem to $m=4l$ for $l \ge 2$?

\section{Final Remarks}
\hspace{4mm}

The matrix case is covered in Subsections 3.2 and 4.2, as the proofs of Theorems \ref{t2} and \ref{t5} need to use Proposition \ref{p0} ($r=1$) and Theorem \ref{t4} ($r=2$), respectively.   This shows that the matrix case ($m=2$) and the higher order tensor case ($m \ge 3$) are connected to each other.

In Subsection 4.2, we allow $r \le 2n-4$ instead of $r \le n$.   Is this a general situation for even $r$?   Further investigation is needed.

From Subsections 3.3 and 4.3, we see that Theorem \ref{t2} may still hold as long as $r$ is odd, even if $GCD(m, r) > 1$; and that Theorem \ref{t5} may still hold as long as $r$ is even, even if $GCD(m, r) > 2$.   Are these true in general?  How can we prove these?   We may see that the proofs of Theorems \ref{t2} and \ref{t5} rely on Theorem \ref{t1}, but the proofs of Theorems \ref{t3} and \ref{t6} do not use a unified technique like Theorem 1.   Can we have a unified technique to study the case that $r$ is odd, $GCD(m, r) > 1$, and the case that $r$ is even, $GCD(m, r) > 2$?

These remain for further investigation.

%\section{The Case that $m=2$ and $n \ge r$: the Matrix Case}
%\hspace{4mm}

%{\bf Guoyin,  It seems that we can get results for $r=3$ and $r=4$.    How can we extend to $r \ge 5$?}

\hspace{2mm}

%%%%%%%%%%%%%%%%%%%%%%%%%%%%%%%%%%%%%%%%%%%%%%%%%%%%%%%%%%%%%%%%%%%%%%%%%%%%%%%%%%%%%%%%%%%%%%%%%%%%%%%%%%%%%%%%%%%%%%%%%%%%%%%%%%%%%%%%%%%%%%%%%%%%%%%%%%%%

\end{document}